# CONCENTRATION OF NORMALIZED SUMS AND A CENTRAL LIMIT THEOREM FOR NONCORRELATED RANDOM VARIABLES[1]

### By Sergey G. Bobkov

### *University of Minnesota*

For noncorrelated random variables, we study a concentration property of the family of distributions of normalized sums formed by sequences of times of a given large length.

**1. Introduction.** Let $X = (X_1, \ldots, X_n)$ be a vector of $n$ random variables on a probability space $(\Omega, \mathbf{P})$ such that, for all $i, j = 1, \ldots, n$,

$$(1.1) \qquad \mathbf{E} X_i X_j = \delta_{ij},$$

where $\delta_{ij}$ is Kronecker's symbol. Given a positive integer $k \leq n$, denote by $\mathcal{G}_{n,k}$ the family of all collections of indices $\tau = \{i_1, \ldots, i_k\}$ of size $k$ with $1 \leq i_1 < \cdots < i_k \leq n$. To every $\tau \in \mathcal{G}_{n,k}$ we associate a normalized sum

$$S_\tau = \frac{X_{i_1} + \cdots + X_{i_k}}{\sqrt{k}}$$

and a corresponding distribution function $F_\tau(x) = \mathbf{P}\{S_\tau \leq x\}$, $x \in \mathbf{R}$. In this paper we show that, when $k$ is a large fixed number, most of the random variables $S_\tau$ are "almost" equidistributed, that is, most of $F_\tau$'s are close to the average distribution function

$$(1.2) \qquad F(x) = \frac{1}{C_n^k} \sum_\tau F_\tau(x),$$

where $C_n^k = \mathrm{card}(\mathcal{G}_{n,k}) = \frac{n!}{k!(n-k)!}$ stands for the usual combinatorial coefficients. To study the rate of closeness, we use the Lévy distance $L(F_\tau, F)$, which is defined to be the infimum over all $\delta \geq 0$ such that $F(x - \delta) - \delta \leq F_\tau(x) \leq F(x + \delta) + \delta$ for all $x \in \mathbf{R}$. In terms of the normalized counting measure $\mu = \mu_{n,k}$ on $\mathcal{G}_{n,k}$, we have:

Received April 2003; revised September 2003.

[1]Supported in part by NSF Grant DMS-01-03929/0405587.

THEOREM 1.1.    *Under* (1.1), *for all* $\delta > 0$,

$$(1.3) \qquad \mu\{\tau : L(F_\tau, F) \geq \delta\} \leq Ck^{3/4}\exp(-ck\delta^8),$$

*where* $C$ *and* $c$ *are certain positive numerical constants.*

The property that, for a growing number of summands $k$, many $F_\tau$'s approximate a "center" $F$ may be viewed as a weak kind of a central limit theorem. In general, however, the center $F$ essentially depends on $k$ and the distribution of the underlying sequence $X$.

An analogous concentration property has been intensively studied in a number of related randomized models. In a seminal work **(year?)**, as an application of the isoperimetric theorem on the sphere, Sudakov established a concentration property of distributions of the weighted sums $\sum_{j=1}^{n} \theta_j X_j$ provided that the weights $\theta_j$ are randomly chosen as coordinates of a point on the unit Euclidean sphere in $\mathbf{R}^n$ (with respect to the uniform measure on the sphere). A different approach in the case of normalized Gaussian weights was suggested by von Weizsäcker **(year?)**. Quantitative versions with refinements for the rate of concentration in the case of log-concave random vectors $X$ were obtained in **(year?)**, **(year?)**; see also **(year?)**, **(year?)**, **(year?)**. Multidimensional random projections of $X$ were considered by Naor and Romik **(year?)**, who essentially used a concentration inequality on the Grassmanian manifold.

As it turns out, the weights can be restricted to the form $\theta_j = \frac{\pm 1}{\sqrt{n}}$ (cf. **(year?)**). As well as on the sphere, the latter model uses a specific dimension-free concentration property on the discrete cube. Similarly, under the conditions of Theorem 1.1 we are dealing with certain weights, namely of the form

$$\theta_j = \frac{\varepsilon_j}{\sqrt{k}}, \qquad 1 \leq j \leq n,$$

where the sequence $(\varepsilon_j)$ contains exactly $k$ 1's and $n - k$ 0's. With respect to the previous examples, this model seems to be closest to the classical, nonrandomized version of the central limit theorem, since only usual sums of data $X_j$ are taken into consideration. The concentration property (1.3) thus tells us that the resulting sum does not depend, in essence, on the concrete times when the observations are made.

Moreover, under certain natural assumptions on random variables $X_j$, the average distribution $F$ must be close to the standard normal distribution function $\Phi$. Namely, suppose we have an infinite sequence of random variables $X_j$ that satisfy the orthogonality condition (1.1).

THEOREM 1.2.    *Let* $\mathbf{E}X_j = 0$ *and* $\sup_j \mathbf{E}|X_j|^3 < \infty$. *Suppose that in probability, as* $n \to \infty$,

$$\frac{X_1^2 + \cdots + X_n^2}{n} \to 1.$$

*Then for all $(k,n)$ such that $1 \ll k \ll n$, for every $\delta > 0$ and for all $\tau \in \mathcal{G}_{n,k}$ except for a set of $\mu$-measure at most $Ck^{3/4}\exp(-ck\delta^8)$, we have $L(F_\tau, \Phi) < \delta + o(1)$.*

Here $o(1)$ denotes a certain sequence $\varepsilon_{n,k}$, independent of $\delta$, which converges to zero for the indicated range of $(k,n)$.

The proofs of Theorems 1.1 and 1.2 are given in Sections 3 and 5. The proof of Theorem 1.1 relies on a concentration property of the measure $\mu$ with respect to the canonical graph structure on $\mathcal{G}_{n,k}$. We discuss this property separately in Section 2. Section 4 is devoted to one auxiliary inequality on elementary symmetric polynomials that is needed for Theorem 1.2. It is also applied in Section 6 to study the asymptotic normality of normalized sums for finite exchangeable sequences.

## 2. Concentration on slices of the discrete cube.

In this section it is convenient to identify $\mathcal{G}_{n,k}$ with the subset of the discrete cube, so let us redefine it as

$$\mathcal{G}_{n,k} = \{x = (x_1, \ldots, x_n) \in \{0,1\}^n : x_1 + \cdots + x_n = k\}.$$

From the discrete cube, $\mathcal{G}_{n,k}$ inherits the structure of a graph: Neighbors are couples of the points which differ exactly in two coordinates. We equip $\mathcal{G}_{n,k}$ with the metric

$$\rho(x,y) = \tfrac{1}{2}\,\mathrm{card}\{i \le n : x_i \ne y_i\}, \qquad x, y \in \mathcal{G}_{n,k},$$

which is one half of the Hamming distance. Every point $x \in \mathcal{G}_{n,k}$ has $k(n-k)$ neighbors $\{s_{ij}x\}_{i \in I(x), j \in J(x)}$ parametrized by

$$I(x) = \{i \le n : x_i = 1\}, \qquad J(x) = \{j \le n : x_j = 0\}.$$

Namely, $(s_{ij}x)_r = x_r$ for $r \ne i, j$ and $(s_{ij}x)_i = x_j$, $(s_{ij}x)_j = x_i$.

For every function $f$ on $\mathcal{G}_{n,k}$ and a point $x$ in $\mathcal{G}_{n,k}$, the discrete gradient $\nabla f(x)$ represents a vector in the Euclidean space $\mathbf{R}^{I(x)} \times \mathbf{R}^{J(x)}$ of dimension $k(n-k)$ with coordinates $(f(x) - f(s_{ij}x))_{i \in I(x), j \in J(x)}$. It has Euclidean length $|\nabla f(x)|$ given by

$$|\nabla f(x)|^2 = \sum_{\rho(x,y)=1} |f(x) - f(y)|^2 = \sum_{i \in I(x)} \sum_{j \in J(x)} |f(x) - f(s_{ij}x)|^2.$$

In 1987, Diaconis and Shahshahani (year?), using a group representation approach, derived a remarkable inequality of Poincaré-type on this graph:

$$(2.1) \qquad \int f^2 \, d\mu - \left( \int f \, d\mu \right)^2 \le \frac{1}{2n} \int |\nabla f|^2 \, d\mu.$$

Note that the constant on the right-hand side can be chosen independently of $k$. Actually, for the quadratic form $(Qf, f) = \int |\nabla f|^2 \, d\mu$ in $L^2(\mathcal{G}_{n,k}, \mu)$, all eigenfunctions and eigenvalues are known. As emphasized in **(year?)**, first they were essentially determined without using group theory by Karlin and McGregor **(year?)**. In particular, with our notations (2.1) becomes equality for all linear functions $f(x) = a_1 x_1 + \cdots + a_n x_n$.

If $|\nabla f|$ is bounded by a constant, say, $\sigma$ (such functions may be viewed as Lipschitz with Lipschitz seminorm at most $\sigma$), then by (2.1), $\mathrm{Var}_\mu(f) \leq \frac{1}{2n} \sigma^2$. This already shows that Lipschitz functions are strongly concentrated around their $\mu$ means $\mathbf{E}_\mu f \equiv \int f \, d\mu$. Applying (2.1) to functions of the form $e^{tf}$ and properly iterating over small $t$, we arrive at a much better estimate,

$$(2.2) \qquad \mu\{|f - \mathbf{E}_\mu f| \geq h\} \leq C e^{-c\sqrt{n}h/\sigma}, \qquad h > 0,$$

up to some numerical positive constants $C$ and $c$. The property that Poincaré-type inequalities imply exponential bounds on the tails of Lipschitz functions was first observed by Gromov and Milman **(year?)** (in the context of Riemannian manifolds) and by Borovkov and Utev **(year?)** (for probability measures on the real line). Afterward it was intensively studied in the literature; see **(year?)** for an extension to the graph setting or **(year?)** for an account of the question.

Although it is not possible to sharpen (2.2) on the basis of (2.1), we may wonder, in analogy with the usual discrete cube, whether a stronger Gaussian bound such as

$$(2.3) \qquad \mu\{|f - \mathbf{E}_\mu f| \geq h\} \leq C \exp(-cnh^2/\sigma^2), \qquad h > 0,$$

holds in the case of the graph $\mathcal{G}_{n,k}$. As is well known, in general, such an improvement can be reached by virtue of a logarithmic Sobolev inequality. An important step in this direction was made by Lee and Yau **(year?)**. They proved that, for every real-valued function $f$ on $\mathcal{G}_{n,k}$,

$$(2.4) \qquad \mathrm{Ent}_\mu(f^2) \leq \frac{C \log(n/k)}{n} \int |\nabla f|^2 \, d\mu,$$

where $C$ is a numerical constant and where we assume for simplicity of notations that $k \leq \frac{n}{2}$. (A little weaker inequality with factor $\log n$ in the place of $\log \frac{n}{k}$ was earlier obtained in **(year?)**.) Here and elsewhere, the entropy functional is defined by

$$\mathrm{Ent}(g) = \mathbf{E}g \log g - \mathbf{E}g \log \mathbf{E}g, \qquad g \geq 0.$$

Thus, when $k$ is proportional to $n$, say, of order $\frac{n}{2}$, the additional logarithmic term $\log \frac{n}{k}$ vanishes and then the logarithmic Sobolev inequality (2.4) represents an improvement, up to a factor in the constant, of the spectral gap inequality (2.1) and implies, in particular, the Gaussian deviation inequality (2.3).

As for the range $k = o(n)$, we have to keep in mind that the constant on the right-hand side of (2.4) is asymptotically sharp. Therefore, to reach (2.3) for the whole range, we need a different argument, and it appears that a modified form of (2.4) may still be used:

THEOREM 2.1. *For every real-valued function $f$ on $\mathcal{G}_{n,k}$,*

$$(2.5) \qquad (n+2)\operatorname{Ent}_\mu(e^f) \leq \mathcal{E}(e^f, f) \leq \int |\nabla f|^2 e^f \, d\mu.$$

*In particular, if $|\nabla f| \leq \sigma$,*

$$(2.6) \qquad \mu\{|f - \mathbf{E}_\mu f| \geq h\} \leq 2\exp(-(n+2)h^2/(4\sigma^2)), \qquad h > 0.$$

The Dirichlet form that appears in the middle of (2.5) is defined canonically by

$$\mathcal{E}(f,g) = \int \langle \nabla f(x), \nabla g(x) \rangle \, d\mu(x)$$

$$= \int \sum_{\rho(x,y)=1} (f(x) - f(y))(g(x) - g(y)) \, d\mu(x),$$

where $f$ and $g$ are arbitrary functions on $\mathcal{G}_{n,k}$. The estimate (2.6) is obtained from (2.5) by applying the latter to functions $tf$: It then yields a distributional inequality $(n+2)\operatorname{Ent}_\mu(e^{tf}) \leq \sigma^2 t^2 \mathbf{E}_\mu e^{tf}$, which is known to imply the bound

$$\mathbf{E}_\mu \exp(t(f - \mathbf{E}_\mu f)) \leq \exp(\sigma^2 t^2/(n+2)), \qquad t \in \mathbf{R},$$

on the Laplace transform of $f$ (an argument due to Ledoux **(year?)**).

The second inequality in (2.5) holds true for the uniform probability measure on an arbitrary finite undirected graph, due to the elementary estimate $(a - b)(e^a - e^b) \leq (a - b)^2(e^a + e^b)/2$, $a, b \in \mathbf{R}$. As for the first inequality in (2.5), it comes naturally in the Markov chain setting in connection with the problem on the rate of convergence to the stationary distribution. In the case of $\mathcal{G}_{n,k}$, it was recently proved **(year?)** in a little more general form by interpolating between the Poincaré and the modified log-Sobolev inequality, and independently **(year?)** where a martingale approach was used to get an asymptotically equivalent constant on the left-hand side of (2.5). For more details and discussions of that inequality, we also refer the interested reader to **(year?)**. Here, for the sake of completeness and to emphasize the "concentration" content, we include below a direct inductive argument.

PROOF OF THEOREM 2.1. For $1 \leq k \leq n-1$, let $A_{n,k}$ denote the best constant in

$$(2.7) \qquad \operatorname{Ent}_\mu(f) \leq A_{n,k}\mathcal{E}(f, \log f) = \frac{A_{n,k}}{C_n^k} \sum_{\rho(x,y)=1} R(f(x), f(y)),$$

where $f$ is an arbitrary positive function on $\mathcal{G} = \mathcal{G}_{n,k}$, $R(a,b) = (a-b)(\log a - \log b)$, for $a, b > 0$ and the summation is performed over all ordered pairs $(x, y) \in \mathcal{G} \times \mathcal{G}$ such that $\rho(x, y) = 1$. By symmetry, $A_{n,k} = A_{n,n-k}$.

When $k = 1$, $\mathcal{G}$ represents a graph of size $n$ where all different points are neighbors of each other (a complete graph). In this case, by Jensen's inequality,

$$\mathrm{Ent}_{\mu}(f) \leq \mathrm{cov}_{\mu}(f, \log f) = \frac{1}{2n^2} \sum_{x \neq y} R(f(x), f(y)) = \frac{1}{2n} \mathcal{E}(f, \log f).$$

Hence, $A_{n,1} \leq \frac{1}{2n}$. As for $k \geq 2$, we deduce a recursive inequality that relates $A_{n,k}$ to $A_{n-1,k-1}$ and then we may proceed by induction. Thus, fix $k \geq 2$ and a positive function $f$ on $\mathcal{G}$ with $\int f \, d\mu = 1$ [this can be assumed in view of the homogeneity of (2.7)]. Introduce subgraphs

$$\mathcal{G}_i = \{x \in \mathcal{G} : x_i = 1\}, \qquad 1 \leq i \leq n,$$

and equip them with uniform probability measures $\mu_i$. Since all $\mathcal{G}_i$ can be identified with $\mathcal{G}_{n-1,k-1}$, we may write the definition (2.7) for these graphs:

$$\int_{\mathcal{G}_i} f \log f \, d\mu_i$$
$$\leq \int_{\mathcal{G}_i} f \, d\mu_i \log \int_{\mathcal{G}_i} f \, d\mu_i + \frac{A_{n-1,k-1}}{C_{n-1}^{k-1}} \sum_{x \in \mathcal{G}_i} \sum_{y \in \mathcal{G}_i, \rho(x,y)=1} R(f(x), f(y)).$$

Put $a_i = \int f \, d\mu_i$. Summing the above inequalities over all $i \leq n$ with weight $\frac{1}{n}$ and making use of $\frac{1}{n} \sum_{i=1}^{n} \mu_i = \mu$, we get

$$(2.8) \quad \begin{aligned} &\int f \log f \, d\mu \\ &\leq \frac{1}{n} \sum_{i=1}^{n} a_i \log a_i + \frac{A_{n-1,k-1}}{nC_{n-1}^{k-1}} \sum_{i=1}^{n} \sum_{x \in \mathcal{G}_i} \sum_{y \in \mathcal{G}_i, \rho(x,y)=1} R(f(x), f(y)). \end{aligned}$$

Since $\frac{1}{n} \sum_{i=1}^{n} a_i = \int f \, d\mu = 1$, the first term in (2.8) is estimated from above, according to the case $k = 1$ in (2.7), by $(A_{n,1}/C_n^1) \sum_{i \neq j} R(a_i, a_j)$. Hence, (2.8) implies

$$\mathrm{Ent}_{\mu}(f) \leq \frac{A_{n,1}}{n} \sum_{i \neq j} R(a_i, a_j) + \frac{A_{n-1,k-1}}{nC_{n-1}^{k-1}} \sum_{i=1}^{n} \sum_{x \in \mathcal{G}_i} \sum_{y \in \mathcal{G}_i, \rho(x,y)=1} R(f(x), f(y)).$$

Now, given $x, y \in \mathcal{G}$ with $\rho(x, y) = 1$, the number of all $i$ such that $x \in \mathcal{G}_i$ and $y \in \mathcal{G}_i$ simultaneously is equal to $k - 1$. Hence, the triple sum contributes

$$(k-1) \sum_{x \in \mathcal{G}} \sum_{y \in \mathcal{G}, \rho(x,y)=1} R(f(x), f(y)) = (k-1)C_n^k \mathcal{E}(f, \log f).$$

Since $((k-1)C_n^k)/(nC_{n-1}^{k-1}) = \frac{k-1}{k}$, we thus get

$$(2.9) \qquad \mathrm{Ent}_\mu(f) \leq \frac{A_{n,1}}{n} \sum_{i \neq j} R(a_i, a_j) + \frac{(k-1)\,A_{n-1,k-1}}{k} \mathcal{E}(f, \log f).$$

To treat the sum in (2.9), note that, for each couple $(i,j)$, $i \neq j$, the map $s_{ij} \colon \{0,1\}^n \to \{0,1\}^n$ acts as a bijection between $\mathcal{G}_i$ and $\mathcal{G}_j$, pushing $\mu_i$ forward onto $\mu_j$ (whenever $k \geq 2$). In particular, $a_j = \int f(y)\,d\mu_j(y) = \int f(s_{ij}x)\,d\mu_i(x)$. Hence, by convexity of $R$ in the positive quarter $a,b > 0$ and Jensen's inequality,

$$R(a_i, a_j) = R\Big( \int f(x)\,d\mu_i(x), \int f(s_{ij}x)\,d\mu_i(x) \Big)$$
$$\leq \int R(f(x), f(s_{ij}x))\,d\mu_i(x).$$

Therefore,

$$(2.10) \qquad \sum_{i \neq j} R(a_i, a_j) \leq \frac{1}{C_{n-1}^{k-1}} \sum_{i \neq j} \sum_{x \in \mathcal{G}_i} R(f(x), f(s_{ij}x)).$$

Note that $y = s_{ij}x$ always implies $\rho(x,y) \leq 1$ and in the case $x \in \mathcal{G}_i$, the equality $\rho(x,y) = 1$ is only possible when $x_i = 1$, $x_j = 0$. Hence, the double sum in (2.10) contains only terms $R(f(x), f(y))$ with $\rho(x,y) = 1$ [the cases $\rho(x,y) = 0$ can be excluded]. In turn, for any couple $x,y \in \mathcal{G}$ such that $\rho(x,y) = 1$, there is a unique pair $(i,j)$ such that $i \neq j$ and $y = s_{ij}x$. Thus, the right-hand side of (2.10) turns into

$$\frac{1}{C_{n-1}^{k-1}} \sum_{x \in \mathcal{G}} \sum_{y \in \mathcal{G}, \rho(x,y)=1} R(f(x), f(y)) = \frac{n}{k} \mathcal{E}(f, \log f)$$

and we finally get, from (2.9),

$$\mathrm{Ent}_\mu(f) \leq \frac{A_{n,1} + (k-1)A_{n-1,k-1}}{k} \mathcal{E}(f, \log f).$$

Hence, $A_{n,k} \leq \frac{1}{k}(A_{n,1} + (k-1)A_{n-1,k-1})$, or $B_{n,k} \leq A_{n,1} + B_{n-1,k-1}$ in terms of $B_{n,k} = kA_{n,k}$. Applying this inequality successively $k-1$ times and recalling that $A_{r,1} \leq \frac{1}{2r}$, we arrive at

$$B_{n,k} \leq \frac{1}{2n} + \frac{1}{2(n-1)} + \cdots + \frac{1}{2(n-(k-2))} + \frac{1}{2(n-(k-1))}.$$

If $k \leq \frac{n}{2}$, each of the above $k$ terms does not exceed $\frac{1}{n+2}$, so $B_{n,k} \leq \frac{k}{n+2}$. This yields the desired estimate $A_{n,k} \leq \frac{1}{n+2}$. In the case $k \geq \frac{n}{2}$, we have $A_{n,k} = A_{n,n-k}$, and Theorem 2.1 follows. $\quad \square$

**3. Proving Theorem 1.1.** We turn to the proof of Theorem 1.1 and to the original definition of $\mathcal{G}_{n,k}$ as a collection of all subsets of $\{1, \ldots, n\}$ of cardinality $k$. We always assume the basic orthonormal hypothesis (1.1) on the sequence $X_1, \ldots, X_n$.

First we focus on the concentration property of the family $\{F_\tau\}$ in terms of their characteristic functions

$$f_\tau(t) = \mathbf{E} e^{itS_\tau}, \qquad \tau \in \mathcal{G}_{n,k}, t \in \mathbf{R},$$

viewed as complex-valued functions on $\mathcal{G}_{n,k}$ with parameter $t$. As a second step, concentration of values of $f_\tau(t)$ around its $\mu$ mean,

$$f(t) = \int f_\tau(t) \, d\mu(\tau) = \int_{-\infty}^{+\infty} e^{itx} \, dF(x),$$

is converted, with the help of standard facts from Fourier analysis, into the concentration property of distributions in the form (1.3). Although this route is different than that in **(year?)** or **(year?)** for the case of the sphere, it has proved to work well on the discrete cube **(year?)** (see also **(year?)**).

LEMMA 3.1. *For every* $t \in \mathbf{R}$, *the function* $\tau \to f_\tau(t)$ *has gradient on* $\mathcal{G}_{n,k}$ *satisfying*

$$|\nabla f_\tau(t)| \le (|t| + t^2) \sqrt{\frac{n}{k}}, \qquad \tau \in \mathcal{G}_{n,k}.$$

PROOF. Every $\tau$ in $\mathcal{G}_{n,k}$ has $k(n-k)$ neighbors in $\mathcal{G}_{n,k}$,

$$\tau_{u,v} = (\tau \setminus \{u\}) \cup \{v\}, \qquad u \in \tau, v \notin \tau.$$

Hence $S_\tau - S_{\tau_{u,v}} = (X_u - X_v)/\sqrt{k}$ and

$$f_\tau(t) - f_{\tau_{u,v}}(t) = \mathbf{E} \exp(itS_\tau)(1 - \exp(-it(X_u - X_v)/\sqrt{k})).$$

Given a complex-valued function $g$ on $\mathcal{G}_{n,k}$, we apply the equivalent representation for the modulus of gradient,

$$|\nabla g(\tau)| = \sup \left| \sum_{u \in \tau} \sum_{v \notin \tau} a_{u,v}(g(\tau) - g(\tau_{u,v})) \right|,$$

where the supremum runs over all collections of complex numbers $a_{u,v}$ such that $\sum_{u \in \tau} \sum_{v \notin \tau} |a_{u,v}|^2 = 1$. In particular, for $g(\tau) = f_\tau(t)$ we have

$$|\nabla f_\varepsilon(t)| = \sup \left| \mathbf{E} \exp(itS_\tau) \sum_{u \in \tau} \sum_{v \notin \tau} a_{u,v}(1 - \exp(-it(X_u - X_v)/\sqrt{k})) \right|$$

$$\le \sup \mathbf{E} \left| \sum_{u \in \tau} \sum_{v \notin \tau} a_{u,v}(1 - \exp(-it(X_u - X_v)/\sqrt{k})) \right|.$$

Using the estimate $|e^{i\alpha} - 1 - i\alpha| \le \frac{1}{2}\alpha^2$ ($\alpha \in \mathbf{R}$), the assumption $\mathbf{E}(X_u - X_v)^2 = 2$ and the identity $\sup \sum_{u \in \tau} \sum_{v \notin \tau} |a_{u,v}| = \sqrt{k(n-k)}$, we can continue to get

$$|\nabla f_\tau(t)| \le \frac{|t|}{\sqrt{k}} \sup \mathbf{E} \left| \sum_{u \in \tau} \sum_{v \notin \tau} a_{u,v}(X_u - X_v) \right|$$

$$+ \frac{t^2}{2k} \sup \mathbf{E} \sum_{u \in \tau} \sum_{v \notin \tau} |a_{u,v}|(X_u - X_v)^2$$

$$= \frac{|t|}{\sqrt{k}} \sup \mathbf{E} \left| \sum_{u \in \tau} \sum_{v \notin \tau} a_{u,v}(X_u - X_v) \right| + \frac{t^2}{k}\sqrt{k(n-k)}.$$

To treat the last double sum, introduce $b_u = \sum_{v \notin \tau} a_{u,v}$ and $c_v = \sum_{u \in \tau} a_{u,v}$, so we can write

$$\sum_{u \in \tau} \sum_{v \notin \tau} a_{u,v}(X_u - X_v) = \sum_{u \in \tau} b_u X_u - \sum_{v \notin \tau} c_v X_v.$$

By (1.1),

$$\mathbf{E} \left| \sum_{u \in \tau} \sum_{v \notin \tau} a_{u,v}(X_u - X_v) \right|^2 = \mathbf{E} \left| \sum_{u \in \tau} b_u X_u \right|^2 + \mathbf{E} \left| \sum_{v \notin \tau} c_v X_v \right|^2$$

$$= \sum_{u \in \tau} |b_u|^2 + \sum_{v \notin \tau} |c_v|^2,$$

but by Cauchy's inequality,

$$|b_u|^2 \le (n-k) \sum_{v \notin \tau} |a_{u,v}|^2, \qquad |c_v|^2 \le k \sum_{u \in \tau} |a_{u,v}|^2,$$

so $\sum_{u \in \tau} |b_u|^2 + \sum_{v \notin \tau} |c_v|^2 \le (n-k) + k = n$. Therefore, once more by Cauchy's inequality,

$$\mathbf{E} \left| \sum_{u \in \tau} \sum_{v \notin \tau} a_{u,v}(X_u - X_v) \right| \le \sqrt{n}.$$

Thus, we arrive at the bound $|\nabla f_\varepsilon(t)| \le |t|\sqrt{\frac{n}{k}} + t^2\sqrt{\frac{n-k}{k}}$, which finishes the proof. $\qquad\square$

COROLLARY 3.2. *For every $t > 0$ and $h > 0$,*

$$(3.1) \qquad \mu \left\{ \tau : \frac{|f_\tau(t) - f(t)|}{t} \ge h \right\} \le 4 \exp\left( \frac{-kh^4}{8(2+h)^2} \right).$$

PROOF.   Indeed, if $t > \frac{2}{h}$, the probablity on the left-hand side is zero, since $|f_\tau(t) - f(t)| \le 2$. In the other case $t \le \frac{2}{h}$, consider the function $g(\tau) = (f_\tau(t) - f(t))/t$. It has $\mu$-mean zero, and according to Lemma 3.1, its modulus of gradient is bounded by $(1 + t)\sqrt{\frac{n}{k}} \le (1 + \frac{2}{h})\sqrt{\frac{n}{k}}$. The same is true for real and imaginary parts $g_1 = \operatorname{Re} g$ and $g_2 = \operatorname{Im} g$. Thus, we are in a position to apply Theorem 2.1 which gives [replacing $n + 2$ with $n$ in (2.6)]

$$\mu\{|g| \ge h\} = \mu\{|g_1|^2 + |g_2|^2 \ge h^2\}$$
$$\le \mu\{|g_1| \ge h/\sqrt{2}\} + \mu\{|g_2| \ge h/\sqrt{2}\}$$
$$\le 4\exp(-kh^2/8(1 + 2/h)^2).$$

Corollary 3.2 follows.   □

By continuity, inequality (3.1) continues to hold in the limit case $t = 0$. Since

$$\mathbf{E}_\mu S_\tau = \int S_\tau \, d\mu(\tau) = \frac{1}{C_n^k} \sum_{i_1 < \cdots < i_k} \frac{X_{i_1} + \cdots + X_{i_k}}{\sqrt{k}} = \sqrt{k}\, \overline{X},$$

where $\overline{X} = (X_1 + \cdots + X_n)/n$, the limiting case becomes

$$(3.2) \qquad \mu\{\tau : |\mathbf{E} S_\tau - \sqrt{k}\mathbf{E}\overline{X}| \ge h\} \le 4\exp(-kh^4/8(2 + h)^2).$$

Thus, under (1.1), the function $g(\tau) = \mathbf{E} S_\tau$ on $\mathcal{G}_{n,k}$ is strongly concentrated around its mean $\mathbf{E}_\mu g = \sqrt{k}\mathbf{E}\overline{X}$.

For the next step, it is important to sharpen inequality (3.1) by making it uniform with respect to the parameter $t$. In other words, we need to control $\sup_{t>0}(|f_\tau(t) - f(t)|)/t$. This can be achieved at the expense of a small deterioration of the bound on the right-hand side of (3.1). Indeed, let us apply (3.1) to points $t_r = rh^2$, $r = 0, 1, \ldots, N = [\frac{2}{h^3}] + 1$, where $[\cdot]$ stands for the integer part of a real number and where the case $r = 0$ is understood as the inequality (3.2). Then we get

$$\mu\left\{\max_{0 \le r \le N} \frac{|f_\tau(t_r) - f(t_r)|}{t_r} \ge h\right\} \le \sum_{r=0}^{N} \mu\left\{\frac{|f_\tau(t_r) - f(t_r)|}{t_r} \ge h\right\}$$
$$(3.3) \qquad\qquad\qquad \le 4(N + 1)\exp\left(\frac{-kh^4}{8(2 + h)^2}\right)$$
$$\le 4\left(\frac{2}{h^3} + 2\right)\exp\left(\frac{-kh^4}{8(2 + h)^2}\right).$$

To involve all remaining values of $t > 0$ in the maximum on the left-hand side, we may assume, as in the proof of Corollary 3.2, that $0 < t \le \frac{2}{h}$. Let $\mathcal{G}(h)$ denote the collection of all $\tau \in \mathcal{G}_{n,k}$ such that $|f_\tau(t_r) - f(t_r)|/t_r < h$ for all $r = 0, 1, \ldots, N$ simultaneously. Recall that $\mathbf{E} S_\tau^2 = 1$, so $|f_\tau'(t)| \le 1$ and $|f_\tau''(t)| \le 1$, and similarly for $f$.

CASE 1.  $0 < t \leq h$. By Taylor's expansion,

$$\frac{f_\tau(t) - f(t)}{t} = i\mathbf{E}(S_\tau - \sqrt{k}\,\overline{X}) + t\int_0^1 (1-v)(f_\tau''(tv) - f''(tv))\,dv.$$

Hence, if $\tau \in \mathcal{G}(h)$ and in particular $|\mathbf{E}\left(S_\tau - \sqrt{k}\,\overline{X}\right)| < h$, we get

$$\frac{|f_\tau(t) - f(t)|}{t} \leq h + t \leq 2h.$$

CASE 2.  $h \leq t \leq \frac{2}{h}$. Pick an index $r = 0, \ldots, N-1$ such that $t_r < t \leq t_{r+1}$. Recalling that $t_{r+1} - t_r = h^2$ and applying the Lipschitz property of $f_\tau$ and $f$, we may write

$$|f_\tau(t) - f(t)| \leq |f_\tau(t) - f_\tau(t_r)| + |f_\tau(t_r) - f(t_r)| + |f(t_r) - f(t)|$$

$$< 2|t - t_r| + t_r h \leq 2h^2 + t_r h < 2h^2 + th \leq 3th.$$

The assumption $t \geq h$ was used on the last step.

Thus, in both cases we obtain that $\tau \in \mathcal{G}(h)$ implies $\sup_{t>0}(|f_\tau(t) - f(t)|/t) < 3h$. Consequently, by (3.3),

$$(3.4) \quad \mu\left\{\sup_{t>0}\frac{|f_\tau(t) - f(t)|}{t} \geq 3h\right\} \leq 4\left(\frac{2}{h^3} + 2\right)\exp\left(\frac{-kh^4}{8(2+h)^2}\right), \qquad h > 0.$$

This is a desired sharpening of (3.1).

PROOF OF THEOREM 1.1.  We use the following observation due to Bohman (year?). Given characteristic functions $\varphi_1$ and $\varphi_2$ of the distribution functions $F_1$ and $F_2$, respectively, if $|\varphi_1(t) - \varphi_2(t)| \leq \lambda t$ for all $t > 0$, then, for all $x \in \mathbf{R}$ and $a > 0$,

$$F_1(x-a) - \frac{2\lambda}{a} \leq F_2(x) \leq F_1(x+a) + \frac{2\lambda}{a}.$$

The particular case $a = \sqrt{2\lambda}$ gives an important relationship,

$$(3.5) \qquad \frac{1}{2}L(F_1, F_2)^2 \leq \sup_{t>0}\frac{|\varphi_1(t) - \varphi_2(t)|}{t},$$

between characteristic functions and the Lévy distance. Therefore, by (3.4) and (3.5),

$$\mu\left\{\frac{1}{2}L(F_\tau, F)^2 \geq 3h\right\} \leq 4\left(\frac{2}{h^3} + 2\right)\exp\left(\frac{-kh^4}{8(2+h)^2}\right).$$

Replacing $6h$ with $\delta^2$ and noticing that only $0 < \delta \leq 1$ should be taken into consideration, we arrive at the estimate

$$\mu\{L(F_\tau, F) \geq \delta\} \leq \frac{C}{\delta^6}\exp(-ck\delta^8), \qquad \delta > 0,$$

with some positive numerical constants $C$ and $c$. On the other hand, in the latter inequality, we may restrict ourselves to values $\delta > c_1 k^{-1/8}$, which make the bound $(C/\delta^6)\exp(-ck\delta^8)$ smaller than 1, and then we arrive at the required inequality (1.3). Theorem 1.1 has been proved. $\square$

**4. Elementary symmetric polynomials.** We turn to the next natural question regarding approximation of the averarge distribution function $F$. According to the definition (1.2), it has characteristic function

$$(4.1) \qquad f(t) = \frac{1}{C_n^k} \mathbf{E} \sum \exp\left(\frac{it X_{i_1}}{\sqrt{k}}\right) \cdots \exp\left(\frac{it X_{i_k}}{\sqrt{k}}\right), \qquad t \in \mathbf{R},$$

with summation over all increasing sequences $1 \le i_1 < \cdots < i_k \le n$. To better understand possible behavior of such sums, introduce normalized elementary symmetric polynomials in $n$ complex variables of degree $k$:

$$\sigma_k(z) = \frac{1}{C_n^k} \sum_{i_1 < \cdots < i_k} z_{i_1} \cdots z_{i_k}, \qquad z = (z_1, \ldots, z_n) \in \mathbf{C}^n.$$

An account on basic results and some other interesting properties of such polynomials can be found in **(year?)**. For our purposes, it is desirable to relate $\sigma_k$ to arithmetic means

$$\bar{z} = \frac{z_1 + \cdots + z_n}{n}.$$

In this section we derive the following statement of independent interest which seems to be absent from the literature (cf. also **(year?)** for a more general scheme).

PROPOSITION 4.1.   *If $|z_j| \le 1$, $j = 1, \ldots, n$, then for all $1 \le k \le n$,*

$$(4.2) \qquad |\sigma_k(z) - \bar{z}^k| \le 6\frac{k-1}{n-1}.$$

Since $|z_j| \le 1$, both quantities satisfy $|\sigma_k(z)| \le 1$ and $|\bar{z}^k| \le 1$, so $|\sigma_k(z) - \bar{z}^k| \le 2$. We can easily refine this bound by applying the polynomial formula

$$\bar{z}^k = \frac{1}{n^k} \sum_{p_1 + \cdots + p_n = k} \frac{k!}{p_1! \cdots p_n!} z_1^{p_1} \cdots z_n^{p_n}$$

$$= \frac{k! C_n^k}{n^k} \sigma_k(z) + \text{remainder}(z).$$

Then we obtain immediately the estimate

$$|\sigma_k(z) - \bar{z}^k| \le 2\left(1 - \frac{k! C_n^k}{n^k}\right).$$

Here, the right-hand side gets small only in the range $k = o(\sqrt{n})$ in which case it is of order $k^2/n$. The bound of order $\frac{k}{n}$ in (4.2) is asymptotically sharp, but its proof requires more sophisticated arguments.

PROOF OF PROPOSITION 4.1. Let $A_{n,k}$ denote maximum of the left-hand side in (4.2) over all possible vectors $z$ with $|z_j| \leq 1$ for all $1 \leq j \leq n$, and let $B_k$ be an optimal constant in

$$A_{n,k} \leq B_k \frac{k-1}{n-1}, \qquad n \geq k.$$

We need a uniform bound on $B_k$. The case $n = 1$ is trivial, since then $A_{n,1} = 0$. If $n = 2$, by simple algebra,

$$\sigma_2(z) - \bar{z}^2 = -\frac{1}{n(n-1)} \sum_{j=1}^{n} (z_j - \bar{z})^2,$$

so

$$|\sigma_2(z) - \bar{z}^2| \leq \frac{1}{n(n-1)} \sum_{j=1}^{n} |z_j - \bar{z}|^2 \leq \frac{1}{n-1}(1 - |\bar{z}|^2) \leq \frac{1}{n-1}.$$

Hence, $A_{n,2} \leq \frac{1}{n-1}$ and $B_2 \leq 1$. To bound the remaining constants, we deduce recursive inequalities that relate $A_{n,k}$ to $A_{n-1,k-1}$ (and then we can argue as in the proof of Theorem 2.1).

Thus let $n \geq k \geq 3$. With every $z \in \mathbf{C}^n$ we associate $n$ vectors in $\mathbf{C}^{n-1}$,

$$z_{(j)} = (z_1, \ldots, z_{j-1}, z_{j+1}, \ldots, z_n), \qquad 1 \leq j \leq n.$$

In what follows we always assume $|z_j| \leq 1$ for all $1 \leq j \leq n$. Let us mention several simple immediate properties and identities:

1. For all $j \leq n$, $|\bar{z}_{(j)}| \leq 1$.
2. We have $\bar{z}_{(j)} - \bar{z} = -(z_j - \bar{z})/(n-1)$.
3. On the other hand, $\bar{z}_{(j)} - \bar{z} = (\bar{z}_{(j)} - z_j)/n$, so we always have $|\bar{z}_{(j)} - \bar{z}| \leq \frac{2}{n}$.
4. We have $\bar{z} = \frac{1}{n} \sum_{j=1}^{n} \bar{z}_{(j)}$.
5. We have $\sigma_k(z) = \frac{1}{n} \sum_{j=1}^{n} z_j \sigma_{k-1}(z_{(j)})$.

From items 4 and 5 we obtain the representation

$$\sigma_k(z) - \bar{z}^k = \frac{1}{n} \sum_{j=1}^{n} z_j (\sigma_{k-1}(z_{(j)}) - \bar{z}_{(j)}^{k-1}) + \frac{1}{n} \sum_{j=1}^{n} z_j (\bar{z}_{(j)}^{k-1} - \bar{z}^{k-1}).$$

Hence,

$$(4.3) \qquad |\sigma_k(z) - \bar{z}^k| \leq A_{n-1,k-1} + \frac{1}{n} \left| \sum_{j=1}^{n} z_j (\bar{z}_{(j)}^{k-1} - \bar{z}^{k-1}) \right|.$$

Thus, our task is to bound the last term on the right-hand side properly. One natural possibility is to use expansion $\bar{z}_{(j)}^{k-1} - \bar{z}^{k-1} = (\bar{z}_{(j)} - \bar{z}) \sum_{u=0}^{k-2} \bar{z}_{(j)}^u \bar{z}^{k-u-2}$. Then by identity 3, $|\bar{z}_{(j)}^{k-1} - \bar{z}^{k-1}| \leq \frac{2(k-1)}{n}$. Applying this estimate in (4.3), we arrive at

$$(4.4) \qquad\qquad A_{n,k} \leq A_{n-1,k-1} + \frac{2(k-1)}{n}.$$

Successive application of this inequality leads to the rough bound $A_{n,k} = O(k^2/n)$. Nevertheless, (4.4) can be useful for small values of $k$. For example, if $k = 3$, we get

$$A_{n,3} \leq A_{n-1,2} + \frac{4}{n} \leq \frac{1}{n-2} + \frac{4}{n}$$
$$= \frac{1}{n-1}\left(5 - \frac{3n-8}{n(n-2)}\right) < \frac{5}{n-1}, \qquad n \geq 3,$$

so $B_3 < \frac{5}{2}$. Similarly, for $k = 4$, $n \geq 4$, by the previous step,

$$(4.5) \qquad A_{n,4} \leq A_{n-1,3} + \frac{6}{n} \leq \frac{1}{n-3} + \frac{4}{n-1} + \frac{6}{n}$$
$$= \frac{1}{n-1}\left(11 - \frac{4n-18}{n(n-3)}\right) < \frac{12}{n-1},$$

so $B_4 < 4$. Hence, $B_k < 4$ for $k \leq 4$, as stated in (4.2).

Thus, assume $n \geq k \geq 5$. We need a more careful estimate of the right-hand side of (4.3) that is independent on $k$. From Taylor's expansion in the integral form, with integration along a segment on the plane connecting two points $a, a_0 \in \mathbf{C}$, we have a canonical estimate

$$|a^{k-1} - a_0^{k-1} - (k-1)a_0^{k-2}(a - a_0)|$$
$$\leq \frac{(k-1)(k-2)}{2}|a - a_0|^2 \max\{|a_0|^{k-3}, |a|^{k-3}\}.$$

In particular, when $a_0 = \bar{z}$, $a = \bar{z}_{(j)}$, we may write, applying property 3,

$$\bar{z}_{(j)}^{k-1} - \bar{z}^{k-1}$$
$$= (k-1)(\bar{z}_{(j)} - \bar{z})\bar{z}^{k-2} + \theta_j \frac{(k-1)(k-2)}{2}|\bar{z}_{(j)} - \bar{z}|^2\left(|\bar{z}| + \frac{2}{n}\right)^{k-3}$$

for some $|\theta_j| \leq 1$. Hence, by statement 2,

$$\sum_{j=1}^{n} z_j(\bar{z}_{(j)}^{k-1} - \bar{z}^{k-1}) = -\frac{k-1}{n-1}\bar{z}^{k-2}\sum_{j=1}^{n} z_j(z_j - \bar{z})$$
$$+ \frac{(k-1)(k-2)}{2(n-1)^2}\left(|\bar{z}| + \frac{2}{n}\right)^{k-3}\sum_{j=1}^{n}\theta_j|z_j - \bar{z}|^2.$$

However, $\sum_{j=1}^{n} z_j(z_j - \bar{z}) = \sum_{j=1}^{n}(z_j - \bar{z})^2$ is bounded in absolute value by $n(1 - |\bar{z}|^2)$. Therefore,

$$(4.6) \quad \frac{1}{n}\left|\sum_{j=1}^{n} z_j(\bar{z}_{(j)}^{k-1} - \bar{z}^{k-1})\right| \leq \frac{k-1}{n-1}|\bar{z}|^{k-2}(1 - |\bar{z}|^2)$$

$$(4.7) \qquad\qquad\qquad + \frac{(k-1)(k-2)}{2(n-1)^2}\left(|\bar{z}| + \frac{2}{n}\right)^{k-3}(1 - |\bar{z}|^2).$$

To bound the expression in (4.6), note that, given $r > 1$, a function of the form $\psi(b) = b^{r-1}(1-b)$ is maximized in $0 \leq b \leq 1$ at $b = 1 - \frac{1}{r}$ and its maximum $(1 - \frac{1}{r})^r \frac{1}{r-1}$ can be bounded by $\frac{1}{e(r-1)}$. Applying this observation with $b = |\bar{z}|^2$ and $r = \frac{k}{2}$, we conclude that

$$(4.8) \qquad \frac{k-1}{n-1}|\bar{z}|^{k-2}(1 - |\bar{z}|^2) \leq \frac{2}{e}\frac{k-1}{k-2}\frac{1}{n-1} < \frac{1}{n-1},$$

where we used the assumption $k \geq 5$.

Next, to bound the expression in (4.7), consider a function of the form $\psi(b) = (b+\varepsilon)^{r-1}(1-b)$ with $\varepsilon = \frac{2}{n}$, $r > 1 + \varepsilon$. It is maximized in $0 \leq b \leq 1$ at $b = 1 - \frac{1+\varepsilon}{r}$ and its maximum $(1+\varepsilon)^r(1 - \frac{1}{r})^r \frac{1}{r-1}$ can be bounded by $((1+\varepsilon)^r)/(e(r-1))$. In particular, with $b = |\bar{z}|$ and $r = k-2$ this yields

$$\frac{(k-1)(k-2)}{2(n-1)^2}\left(|\bar{z}| + \frac{2}{n}\right)^{k-3}(1 - |\bar{z}|)(1 + |\bar{z}|)$$

$$\leq \frac{(k-1)(k-2)}{e(k-3)(n-1)^2}\left(1 + \frac{2}{n}\right)^{k-2}.$$

Hence, using $(1 + \frac{2}{n})^{k-2} \leq (1 + \frac{2}{n})^{n-2} \leq e^2/(1 + \frac{2}{n})^2$ and $\frac{(k-1)(k-2)}{k-3} = k + \frac{2}{k-3} < n+2$, we can estimate the expression in (4.7) by $\frac{e}{n-1} < \frac{3}{n-1}$. Together with (4.8), the left-hand side of (4.6) is thus bounded by $\frac{4}{n-1}$. Thus, returning to (4.3), we obtain a more precise recursive inequality than (4.4):

$$(4.9) \qquad A_{n,k} \leq A_{n-1,k-1} + \frac{4}{n-1}, \qquad n \geq k \geq 5.$$

Finally, applying (4.9) $k-4$ times and the obtained estimate (4.5), we get

$$(4.10) \quad A_{n,k} \leq \frac{4}{n-1} + \frac{4}{n-2} + \cdots + \frac{4}{n-k+4} + A_{n-k+4,4} \leq \frac{4(k-1)}{n-k+3}.$$

In the case $k \leq \frac{n}{3} + 1$, we have $n - k + 3 \geq \frac{2}{3}n$, and (4.10) yields the desired estimate (4.2). In the other case there is nothing to prove since then $\frac{6(k-1)}{n-1} \geq 2 \geq A_{n,k}$. Proposition 4.1 is proved. $\quad\square$

**5. Theorem 1.2 and its generalization.** As before, let $X_1, \ldots, X_n$ be random variables that satisfy the orthogonality condition (1.1) and let $1 \le k \le n$. Now we are prepared to study asymptotic properties of the average characteristic function $f$ defined in (4.1) (of the average distribution function $F$). Given $\omega \in \Omega$, introduce random characteristic functions

$$f_\omega(t) = \frac{1}{C_n^k} \sum \exp\left(\frac{itX_{i_1}(\omega)}{\sqrt{k}}\right) \cdots \exp\left(\frac{itX_{i_k}(\omega)}{\sqrt{k}}\right),$$

$$g_\omega(t) = \left(\frac{\exp(itX_1(\omega)) + \cdots + \exp(itX_n(\omega))}{n}\right)^k,$$

where summation runs over all increasing sequences $1 \le i_1 < \cdots < i_k \le n$. Thus, $f(t) = \mathbf{E}\, f_\omega(t)$. Also put

$$g(t) = \mathbf{E}\, g_\omega(t) = \mathbf{E}\left(\frac{\exp(itX_1) + \cdots + \exp(itX_n)}{n}\right)^k, \qquad t \in \mathbf{R}.$$

By Proposition 4.1, we always have $|f_\omega(t) - g_\omega(t)| \le \frac{6k}{n}$, so a similar inequality must hold for corresponding means, that is,

$$(5.1) \qquad\qquad |f(t) - g(t)| \le \frac{6k}{n}, \qquad t \in \mathbf{R}.$$

Hence, when $k = o(n)$, the associated distribution functions must be also close to each other and we may concentrate on the asymptotic behavior of $g$, only.

A probabilistic meaning of each function $g_\omega$ is very simple. Indeed, given $\omega \in \Omega$, let $Y_1, \ldots, Y_k$ be independent, identically distributed random variables defined on some probability space $(M, Q)$, whose common distribution is a sample distribution:

$$Q\{Y_1 = X_j(\omega)\} = \frac{1}{n}, \qquad 1 \le j \le n.$$

Then, by the very definition, $g_\omega$ represents the characteristic function of the random variable

$$T_\omega = \frac{Y_1 + \cdots + Y_k}{\sqrt{k}}.$$

It has $Q$ mean $\mathbf{E}_Q T_\omega = \sqrt{k}\,\overline{X}(\omega)$, where $\overline{X}(\omega) = \frac{1}{n}\sum_{j=1}^n X_j(\omega)$ is just a sample mean associated to the "sample" $X_1, \ldots, X_n$, and has $Q$ variance

$$(5.2) \qquad\qquad \sigma^2(\omega) = \mathrm{Var}_Q(T_\omega) = \frac{1}{n}\sum_{j=1}^n (X_j(\omega) - \overline{X}(\omega))^2,$$

representing the usual sample variance. For simplicity, in some places we omit $\omega$, hoping this does not lead to confusion.

By the canonical central limit theorem, the random variable $T_\omega$ has a distribution function, $G_\omega$, which is close to the normal $N(\sqrt{k}\,\overline{X}, \sigma^2)$. Hence, the distribution function $G(x) = \mathbf{E}\,G_\omega(x)$, associated to the characteristic function $g$, is close to a $\mathbf{P}$ mixture of $N(\sqrt{k}\,\overline{X}, \sigma^2)$-distribution functions. Clearly, this mixture can be described as the distribution function of a random variable of the form

$$\xi = \sqrt{k}\,\overline{X} + \sigma\zeta,$$

where $\zeta$ is a standard normal random variable independent of all r.v.'s $X_j$. It has characteristic function

(5.3) $$h(t) = \mathbf{E}\,e^{it\xi} = \mathbf{E}\exp(\sqrt{k}\,\overline{X}it - \sigma^2 t^2/2).$$

LEMMA 5.1. *If* $\mathbf{E}X_j = 0$, $\mathbf{E}|X_j|^3 \le \beta$ $(1 \le j \le n)$, *then*

$$\sup_{t>0} \frac{|f(t) - h(t)|}{t} \le 3\left(\frac{k}{n}\right)^{1/2} + 6\frac{\beta^{1/4}}{k^{1/8}}.$$

Let $H$ denote the distribution function of $\xi$. By Bohman's inequality (3.5) applied to $F_1 = F$ and $F_2 = H$, we get

$$L^2(F, H) \le 6\left(\frac{k}{n}\right)^{1/2} + 12\frac{\beta^{1/4}}{k^{1/8}}.$$

The quantity on the right-hand side is small once $k$ is large and $\frac{k}{n}$ is small. In this case, $\mathbf{E}(\sqrt{k}\,\overline{X})^2 = \frac{k}{n}$ is small, as well, and according to (5.3), $h(t)$ is close to the characteristic function

$$\mathbf{E}\exp\left\{-\frac{t^2}{2n}\sum_{j=1}^n X_j^2\right\}.$$

Thus, we arrive at the following conclusion that includes the statement of Theorem 1.2. Let $(X_j)_{j=1}^\infty$ be a sequence of random variables that satisfies the correlation condition (1.1) and such that $\mathbf{E}X_j = 0$, $\sup_j \mathbf{E}|X_j|^3 < +\infty$. Assume that, for some random variable $R \ge 0$, as $n \to \infty$,

$$\frac{1}{n}\sum_{j=1}^n X_j^2 \to R^2$$

in the sense of the weak convergence of distributions on the real line. Let $\Phi_R$ denote the distribution function of the random variable $R\zeta$, where $\zeta$ is a standard normal random variable that is independent on $R$. Then we have:

THEOREM 5.2. *For all* $(k, n)$ *in the range* $1 \ll k \ll n$, *for every* $\delta > 0$ *and for all* $\tau \in \mathcal{G}_{n,k}$ *except for a set of* $\mu$ *measure at most* $Ck^{3/4}\exp(-ck\delta^8)$, *we have*

$$L(F_\tau, \Phi_R) < \delta + o(1).$$

PROOF OF LEMMA 5.1.   We use the following standard estimate (needed for the Berry–Esseen theorem; cf., e.g., **(year?)**, Chapter V, paragraph 2, Lemma 1): If $Z_1, \ldots, Z_k$ are independent r.v.'s such that $\mathbf{E} Z_l = 0$, $\mathbf{E} |Z_l|^3 < \infty$ and $B = \sum_{l=1}^{k} \mathbf{E} Z_l^2$, then

$$\left| \mathbf{E} \exp\left( \frac{it(Z_1 + \cdots + Z_k)}{\sqrt{B}} \right) - \exp\left( \frac{-t^2}{2} \right) \right| \leq 16 L |t|^3 \exp\left( \frac{-t^2}{3} \right), \qquad |t| \leq \frac{1}{4L},$$

where $L = B^{-3/2} \sum_{l=1}^{k} \mathbf{E} |Z_l|^3$ (the so-called Lyapunov fraction). Dividing by $t$ and maximizing the right-hand side over all $t > 0$, we get

$$(5.4) \qquad \frac{|\mathbf{E} \exp(it(Z_1 + \cdots + Z_k)/\sqrt{B}) - \exp(-t^2/2)|}{t} \leq 18 L,$$

provided that $0 < t \leq \frac{1}{4L}$. In the case $t \geq \frac{1}{4L}$, the left-hand side can be estimated by $\frac{2}{t} \leq 8L$, so (5.4) holds for all $t > 0$. In particular, if the $Z_l$'s are identically distributed with $\mathbf{E} Z_1^2 = \sigma^2$ and $\mathbf{E} |Z_1|^3 = \beta$, then $B = \sigma^2 k$, $L = \beta / \sigma^3 \sqrt{k}$, and the above bound yields

$$\max_{t > 0} \frac{|\mathbf{E} \exp(it(Z_1 + \cdots + Z_k)/\sqrt{k}) - \exp(-\sigma^2 t^2/2)|}{t} \leq \frac{18}{\sqrt{k}} \, \frac{\beta}{\sigma^3}.$$

In particular, this inequality can be applied on the probability space $(M, Q)$ to random variables $Z_l = Y_l - \overline{X}$. In this case, $\sigma^2 = \sigma^2(\omega)$ represents the sample variance (5.2) and similarly $\beta(\omega) = \frac{1}{n} \sum_{j=1}^{n} |X_j - \overline{X}|^3$. Thus, introducing the characteristic function $h_\omega(t) = \exp(\sqrt{k} \overline{X} it - \sigma^2 t^2/2)$, we obtain that

$$(5.5) \qquad \sup_{t > 0} \frac{|g_\omega(t) - h_\omega(t)|}{t} \leq \frac{18}{\sqrt{k}} \frac{\beta(\omega)}{\sigma^3}.$$

Note that both $g_\omega$ and $h_\omega$ correspond to distributions with expectation $\sqrt{k} \overline{X}$ and variance $\sigma^2$. Hence, by Taylor's expansion around zero, $|g_\omega(t) - h_\omega(t)| \leq \sigma^2 t^2$ for all $t \in \mathbf{R}$. On the other hand, we always have a trivial bound $|g_\omega(t) - h_\omega(t)| \leq 2$. Combining these, we get

$$\frac{|g_\omega(t) - h_\omega(t)|}{t} \leq \min\left\{ \sigma^2 t, \frac{2}{t} \right\} \leq \sqrt{2} \sigma, \qquad t > 0.$$

Together with (5.5) and maximizing over $\sigma > 0$, this gives, for all $t > 0$,

$$\frac{|g_\omega(t) - h_\omega(t)|}{t} \leq \min\left\{ \frac{18}{\sqrt{k}} \frac{\beta(\omega)}{\sigma^3}, \sqrt{2} \sigma \right\} \leq \frac{3\beta(\omega)^{1/4}}{k^{1/8}}.$$

Averaging over $\omega$ and using Hölder's inequality, we obtain that

$$\frac{|g(t) - h(t)|}{t} \leq \frac{3}{k^{1/8}} (\mathbf{E} \beta(\omega))^{1/4}, \qquad t > 0,$$

since $h(t) = \mathbf{E}h_\omega(t)$. To estimate $\mathbf{E}\beta(\omega)$, we may apply Jensen's inequality, implying $|X_j - \overline{X}|^3 \leq \frac{1}{n}\sum_{l=1}^{n}|X_j - X_l|^3$. Since $\mathbf{E}|X_j - X_l|^3 \leq 4\mathbf{E}|X_j|^3 + 4\mathbf{E}|X_l|^3 \leq 8\beta$, we arrive at $\mathbf{E}|X_j - \overline{X}|^3 \leq 8\beta$ and, therefore, $\mathbf{E}\beta(\omega) \leq 8\beta$. Hence,

$$(5.6) \qquad \sup_{t>0}\frac{|g(t)-h(t)|}{t} \leq \frac{6\beta^{1/4}}{k^{1/8}}.$$

It remains to involve the characteristic function $f$. Combining (5.1) and (5.6), we get

$$(5.7) \qquad |f(t)-h(t)| \leq \frac{6k}{n} + \frac{6\beta^{1/4}}{k^{1/8}}t, \qquad t>0.$$

On the other hand, $\mathbf{E}\xi = 0$ and, by independence of $\zeta$ and $(X_1, \ldots, X_n)$,

$$\mathbf{E}\xi^2 = \mathbf{E}(\sqrt{k}\,\overline{X} + \sigma\zeta)^2 = k\mathbf{E}(\overline{X})^2 + \mathbf{E}\sigma^2 = 1 + \frac{k-1}{n} \leq 2,$$

so $h'(0) = 0$ and $|h''(t)| \leq 2$ for all $t \in \mathbf{R}$. In addition, the distribution function $F$ has mean 0 and variance 1, so $f'(0) = 0$ and $|f''(t)| \leq 1$. Consequently, by Taylor's expansion around zero, $\frac{|f(t)-h(t)|}{t} \leq \frac{3t}{2}$, $t > 0$. Together with (5.7), the latter gives

$$\frac{|f(t)-h(t)|}{t} \leq \frac{3}{2}\min\left\{t, \frac{4k}{nt} + \frac{4\beta^{1/4}}{k^{1/8}}\right\}.$$

Finally, let us note that, given $a, b > 0$, a function of the form $u(t) = \min\{t\frac{b}{t} + a\}$ attains its maximum at $t_0 = (a + \sqrt{a^2 + 4b})/2$, and, at this point, $u(t_0) = t_0 \leq a + \sqrt{b}$. Applying this to $b = \frac{4k}{n}$ and $a = 4\beta^{1/4}/k^{1/8}$, we arrive at

$$\sup_{t>0}\frac{|f(t)-h(t)|}{t} \leq \frac{3}{2}\left[2\left(\frac{k}{n}\right)^{1/2} + \frac{4\beta^{1/4}}{k^{1/8}}\right].$$

Lemma 5.1 and, therefore, Theorem 5.2 are proved. $\square$

**6. Exchangeable random variables.** Random variables $X_1, \ldots, X_k$ are called exchangeable (or interchangeable) if the distribution $\mathbf{P}_X$ of the random vector $X = (X_1, \ldots, X_k)$, as a measure on $\mathbf{R}^k$, is invariant under permutations of coordinates. A similar definition applies in the case of an infinite sequence $\{X_k\}_{k=1}^{\infty}$. In particular, for all $k \geq 1$, the distributions of the normalized sums $(X_{i_1} + \cdots + X_{i_k})/\sqrt{k}$ do not depend on the choice of indices $i_1 < \cdots < i_k$. So let

$$S_k = \frac{X_1 + \cdots + X_k}{\sqrt{k}}.$$

Given that

$$\mathbf{E}X_1 = 0, \qquad \mathbf{E}X_1^2 = 1, \tag{6.1}$$

a well-known theorem due to Blum, Chernoff, Rosenblatt and Teicher [4] asserts that $S_k \to N(0,1)$ weakly in distribution as $k \to \infty$ if and only if

$$\mathbf{E}X_1X_2 = 0, \qquad \mathbf{E}X_1^2X_2^2 = 1; \tag{6.2}$$

that is, $\operatorname{cov}(X_1, X_2) = \operatorname{cov}(X_1^2, X_2^2) = 0$. Moreover, Berry–Esseen's bound

$$\sup_{x \in \mathbf{R}} |\mathbf{P}\{S_k \leq x\} - \Phi(x)| \leq c\frac{\mathbf{E}|X_1|^3}{\sqrt{k}}, \tag{6.3}$$

with some universal $c$, extends from the i.i.d. case to this case as well.

Weaker assumptions than (6.1) and (6.2) with different normalization of the sums may also lead to asymptotic normality (see, e.g., [22, 24]). However, less seems to be known in the case of finite sequences of exchangeable variables. A basic tool that allows study of the various properties of an infinite exchangeable sequence $X = \{X_k\}_{k=1}^\infty$ is de Finetti's representation of the distribution $\mathbf{P}_X$ of $X$ as a mixture

$$\mathbf{P}_X = \int_\Pi \mu_\alpha^\infty \, d\pi(\alpha) \tag{6.4}$$

of product probability measures $\mu_\alpha^\infty = \mu_\alpha \otimes \mu_\alpha \otimes \cdots$ on $\mathbf{R}^\infty$. Here $(\Pi, \pi)$ is some probability space and $\{\mu_\alpha\}_{\alpha \in \Pi}$ is some family of marginals with the property that functions $\alpha \to \mu_\alpha(B)$ are $\pi$-measurable for all Borel sets $B$ on the real line. In terms of this representation and assuming (6.1) is fulfilled, the central limit theorem $S_k \to N(0,1)$ holds true if and only if the measures $\mu_\alpha$ have mean 0 and variance 1 for $\pi$-almost all $\alpha$ ([4], Lemma 1). The latter is also characterized directly in terms of $X$ in the form (6.2).

In the general case of a finite exchangeable sequence $X = (X_1, \ldots, X_k)$, the finite-dimensional analogue of (6.4),

$$\mathbf{P}_X(B) = \int \mu_\alpha^k(B) \, d\pi(\alpha), \qquad B \subset \mathbf{R}^k, \tag{6.5}$$

where $\mu_\alpha^k = \mu_\alpha \otimes \cdots \otimes \mu_\alpha$ are product probability measures on $\mathbf{R}^k$, is no longer valid and, in fact, the class of distributions on $\mathbf{R}^k$ invariant under permutations of coordinates is much wider. Therefore, it is natural to associate to $X$ a maximum natural number $n = n(X)$ such that, for some exchangeable sequence $\widetilde{X}_1, \ldots, \widetilde{X}_n$ defined perhaps on a different probability space, the random vectors $(X_1, \ldots, X_k)$ and $(\widetilde{X}_1, \ldots, \widetilde{X}_k)$ are equidistributed. If $n$ can be chosen as large as we wish or, equivalently, if $\mathbf{P}_X$ admits representation (6.5), put $n(X) = \infty$.

It may occur that $X$ has no exchangeable extension: $n(X) = k$. In that case, it is hardly possible to reach asymptotic normality of the normalized sum $S_k$, even under moment assumptions such as (6.1) and (6.2). However, when $n(X) \gg k$, the situation changes considerably. In view of de Finetti's theorem, it seems natural to expect in this case that $\mathbf{P}_X$ has to be close in some sense to the class $\mathcal{M}_k$ of mixtures of product probability measures on $\mathbf{R}^k$. That is, there should hold an approximate equality in (6.5). In terms of the variational distance $\| \cdot \|_{\mathrm{TV}}$ between probability measures, this question was studied by Diaconis and Freedman [16]. It was shown, in particular that, for some $Q$ in $\mathcal{M}_k$,

$$(6.6) \qquad \frac{1}{2} \| \mathbf{P}_X - Q \|_{\mathrm{TV}} \leq 1 - \frac{k! C_n^k}{n^k}, \qquad n = n(X),$$

and that the bound cannot be improved. Actually, if an exchangeable extension $X_1, \ldots, X_n$ exists on the same probability space $(\Omega, \mathbf{P})$, we can take $Q(B) = \int \mu_\omega^k(B) \, d\mathbf{P}(\omega)$, that is, with

$$\Pi = \Omega, \pi = \mathbf{P}, \qquad \mu_\omega = \frac{\delta_{X_1(\omega)} + \cdots + \delta_{X_n(\omega)}}{n}.$$

Under the product measures, the distribution of the function $x \to (x_1 + \cdots + x_k)/\sqrt{k}$ is nearly normal (under proper moment conditions), so the inequality (6.6) can be used to study the asymptotic normality of $S_k$. However, as emphasized in [16], the expression on the right-hand side in (6.6) is of order $k^2/n$ for $k = o(\sqrt{n})$, while it is of order 1 for larger values of $k$. Hence, only the range $k = O(\sqrt{n})$ can be taken into consideration or other metrics that better react on the weak convergence of distributions should be examined in the case $k > O(\sqrt{n})$. In part concerning half-spaces of the form $B = \{x \in \mathbf{R}^k : x_1 + \cdots + x_k \leq c\}$, the closeness of $\mathbf{P}_X(B)$ to $Q(B)$ can be estimated by virtue of Proposition 4.1. As a consequence, we can derive:

PROPOSITION 6.1. *Let $X = (X_1, \ldots, X_k)$, $k \geq 2$, be an exchangeable sequence that satisfies the moment hypotheses* (6.1) *and* (6.2). *Then*

$$(6.7) \qquad \sup_{x \in \mathbf{R}} |\mathbf{P}\{S_k \leq x\} - \Phi(x)| \leq c \left[ \left( \frac{k}{n(X)} \right)^p + \frac{(\mathbf{E}|X_1|^4)^{1/6}}{k^q} \right],$$

*for some universal $c > 0$ and $p, q > 0$.*

Although this is not as sharp as (6.3), we can still control closeness to normality for finite sequences under the same hypotheses. The assumption $\mathbf{E} X_1^4 < +\infty$ is technical and can be a little relaxed (to the third moment, e.g.). The second assumption in (6.2) can be weakened to $\mathbf{E} X_1^2 X_2^2 \leq 1$. Although a strict inequality is impossible here for infinite exchangeable sequences, it does hold for some interesting finite exchangeable sequences (cf., e.g., **(year?)**).

PROOF OF PROPOSITION 6.1. Let $X$ have an exchangeable extension $X_1, \ldots, X_n$ on $(\Omega, \mathbf{P})$. By exchangeability, $F(x) = \mathbf{P}\{S_k \leq x\}$ represents the average distribution function (1.1), and its characteristic function $f$ appears in (4.1). Note that, under the measure $Q(B) = \int \mu_\omega^k(B) \, d\mathbf{P}(\omega)$, the function $x \to (x_1 + \cdots + x_k)/\sqrt{k}$ has distribution $G$ considered along the proof of Theorem 5.2. Moreover, by Lemma 5.1 and Hölder's inequality,

$$\sup_{t>0} \frac{|f(t) - h(t)|}{t} \leq 3\left(\frac{k}{n}\right)^{1/2} + 6\frac{(\mathbf{E}X_1^4)^{1/4}}{k^{1/8}},$$

where we recall that $h(t) = \mathbf{E}\exp(\sqrt{kX}it - \sigma^2 t^2/2)$ represents the characteristic function of $\xi = \sqrt{kX} + \sigma\zeta$ with $\zeta \in N(0,1)$ independent of $(X_1, \ldots, X_n)$. By Bohman's inequality (5.3) and using $\sqrt{a+b} \leq \sqrt{a} + \sqrt{b}$ $(a, b \geq 0)$, we may write down a bound on the Lévy distance,

$$(6.8) \qquad L(F, H) \leq \sqrt{6}\left(\frac{k}{n}\right)^{1/4} + \sqrt{12}\frac{(\mathbf{E}X_1^4)^{1/8}}{k^{1/16}}$$

for the associated distribution functions. Note that we have used the assumptions $\mathbf{E}X_1 = \mathbf{E}X_1X_2 = 0$ and $\mathbf{E}X_1^2 = 1$ in this step.

To quantify closeness of the distribution function $H$ to $\Phi$, we write $\xi = \zeta + \eta$ with a small "error" $\eta = \sqrt{k}\,\overline{X} + (\sigma - 1)\zeta$. We apply the following general observation: For all random variables $\zeta$ and $\eta$,

$$(6.9) \qquad L(F_{\zeta+\eta}, F_\zeta) \leq (\mathbf{E}\eta^2)^{1/3},$$

where $F_{\zeta+\eta}$ and $F_\zeta$ are corresponding distribution functions. Indeed, there is nothing to prove if $\delta \equiv (\mathbf{E}\eta^2)^{1/3} \geq 1$. In the other case, since for all $x \in \mathbf{R}$ and $h > 0$,

$$\{\zeta \leq x\} = \{\zeta \leq x, \eta \leq h\} \cup \{\zeta \leq x, \eta > h\} \subset \{\zeta + \eta \leq x + h\} \cup \{\eta > h\},$$

by Chebyshev's inequality, we get $F_\zeta(x) \leq F_{\zeta+\eta}(x+h) + \mathbf{E}\eta^2/h^2$. Applying the latter to another couple of random variables $(\zeta + \eta, -\eta)$ and to $x - h$ in the place of $x$, we also have $F_{\zeta+\eta}(x-h) \leq F_\zeta(x) + \mathbf{E}\eta^2/h^2$. All this together with $h = \delta$ yields

$$F_{\zeta+\eta}(x-\delta) - \delta \leq F_\zeta(x) \leq F_{\zeta+\eta}(x+h) + \delta,$$

which is exactly (6.9).

Thus, returning to our specific random variables $(\zeta, \eta)$, since $F_\xi = H$ and $F_\zeta = \Phi$, we may conclude that

$$(6.10) \qquad L(H, \Phi) \leq (\mathbf{E}\eta^2)^{1/3}.$$

Now, since $\mathbf{E}X_1X_2 = 0$ and $\mathbf{E}X_1^2 = 1$, we have $\mathbf{E}\eta^2 = \mathbf{E}(\sqrt{k}\,\overline{X})^2 + \mathbf{E}(\sigma - 1)^2 = \frac{k-1}{n} + 2\mathbf{E}(1-\sigma)$. Also note $1 - \sigma = (1 - \sigma^2)/(1+\sigma) = (1 - \overline{X}^2 + (\overline{X})^2)/(1+\sigma)$, so $|1 - \sigma| \le |\overline{X}^2 - 1| + (\overline{X})^2$, where $\overline{X}^2 = \frac{1}{n}\sum_{j=1}^{n} X_j^2$. By the assumption $\mathbf{E}X_1^2 X_2^2 = 1$ and since $\mathbf{E}\overline{X}^2 = 1$,

$$\mathbf{E}|\overline{X}^2 - 1|^2 = \mathrm{Var}(\overline{X}^2) = \frac{\mathbf{E}X_1^4}{n} + \frac{n(n-1)\,\mathrm{cov}(X_1^2, X_2^2)}{n} = \frac{\mathbf{E}X_1^4}{n}.$$

Therefore, $\mathbf{E}|\overline{X}^2 - 1| \le (\mathbf{E}X_1^4)^{1/2}/\sqrt{n}$ and $\mathbf{E}|1 - \sigma| \le (\mathbf{E}X_1^4)^{1/2}/\sqrt{n} + \frac{1}{n}$. Thus, we get $\mathbf{E}\eta^2 \le (k+1)/n + (2(\mathbf{E}X_1^4)^{1/2})/\sqrt{n}$ and by (6.10),

$$L(H, \Phi) \le 2\left(\frac{k}{n}\right)^{1/3} + \frac{2(\mathbf{E}X_1^4)^{1/6}}{n^{1/6}}.$$

Combining this with (6.8) and making use of $k \le n$ and $(\mathbf{E}X_1^4)^{1/8} \le (\mathbf{E}X_1^4)^{1/6}$ (since the fourth moment is greater than or equal to 1), we obtain that

$$L(F, \Phi) \le c_1\left(\frac{k}{n}\right)^{1/4} + c_2\frac{(\mathbf{E}X_1^4)^{1/6}}{k^{1/16}}.$$

Finally, we always have $\|F - \Phi\|_\infty \le 2\,L(F, \Phi)$, so we arrive at (6.7) with $p = \frac{1}{4}$, $q = \frac{1}{16}$. This completes the proof. $\square$

SCHOOL OF MATHEMATICS
UNIVERSITY OF MINNESOTA
127 VINCENT HALL
206 CHURCH ST. S.E.
MINNEAPOLIS, MINNESOTA 55455
USA
E-MAIL: bobkov@math.umn.edu